\documentstyle{article}

\newtheorem{theorem}{Theorem}

\def\QED{\quad\blackslug\lower 8.5pt\null}

\begin{document}

\begin{center}
{\large   THE DARBOUX MAPPING OF }
\vspace*{2mm}
{\large  CANAL HYPERSURFACES}\\ 
\vspace*{7mm}
{\small BY} \\
\vspace*{3mm}
 Maks A. Akivis\\
\vspace*{3mm}
{\footnotesize {\em Department of Mathematics and Computer Science\\
Ben-Gurion University of the Negev\\
P. O. Box 653, Beer-Sheva, 84105, Israel\\
e-mail: akivis@black.bgu.ac.il}}\\
\vspace*{5mm}
{\small AND} \\
\vspace*{3mm}
 Vladislav V. Goldberg\\
\vspace*{3mm}
{\footnotesize {\em Department of Mathematics \\
New Jersey Institute of Technology \\
Newark, NJ 07102, USA\\
e-mail: vlgold@numerics.njit.edu}}\\
\end{center}
\vspace*{3mm}
\begin{abstract}
The geometry of canal hypersurfaces of an $n$-dimensional 
conformal space $C^n$ is studied. Such hypersurfaces 
are envelopes of $r$-parameter families of hyperspheres, 
 $1 \leq r \leq n-2$. In the present paper the conditions that characterize 
canal hypersurfaces, and which were known earlier, are made more precise.
The main attention is given to the study of the Darboux maps 
of canal hypersurfaces in the de Sitter space $M_1^{n+1}$ 
and the projective space $P^{n+1}$. To canal hypersurfaces 
there correspond $r$-dimensional spacelike tangentially 
nondegenerate 
submanifolds in $M_1^{n+1}$ and tangentially degenerate 
hypersurfaces of rank $r$ in $P^{n+1}$. In this connection 
the problem of existence of singular points on canal 
hypersurfaces is considered.
\end{abstract}

\vspace*{15mm}

\setcounter{equation}{0}

\setcounter{section}{0}

{\bf 0. Introduction.}  
Canal hypersurfaces are envelopes of families of hyperspheres. 
In a three-dimensional Euclidean space canal surfaces were 
considered in many classical texts on differential geometry 
(see, for example, the book [B 29]). Canal hypersurfaces in a 
multidimensional Euclidean space were also investigated in many 
papers. Since the property of a hypersurface to be a canal 
hypersurface is conformally invariant, it is natural to consider 
the canal hypersurfaces in a conformal space $C^n$. Such 
investigations have been done,  for example, in 
[A 52],  [M 55], and [V 57, 62]. 

However, in all these works the authors did not note the 
difference between 1-canal hypersurfaces (i.e., the envelopes of 
$(n-2)$-parameter families of hyperspheres) and $m$-canal 
hypersurfaces, where $m \geq 2$ (i.e., the envelopes of 
$(n-m-1)$-parameter families of hyperspheres).

In the present paper we discuss this difference, namely the fact 
that the analytic characterization of 1-canal hypersurfaces is 
connected with a third-order differential neighborhood of a point 
$x$ of a hypersurface $V^{n-1}$ while that of  $m$-canal 
hypersurfaces for $m \geq 2$ is connected with a second-order 
differential neighborhood of $x \in V^{n-1}$. 

After recalling the fundamental equations of the conformal theory 
of hypersurfaces (Section {\bf 1}), we find separately an 
analytic characterization of $m$-canal hypersurfaces (Section 
{\bf 2}) and of  1-canal hypersurfaces (Section {\bf 3}) and 
compare them.  Note that while Theorem 1 and Theorem 3 were known 
(see [A 52], [V 57]   and [AG 96], p. 61,  
for Theorem 1 and [A 65]  and [AG 96], p. 58, for Theorem 3), 
Theorem 2 in the form as presented in this paper appears for the 
first time. 

In Sections {\bf 4} and {\bf 5} 
we study the Darboux mapping of canal hypersurfaces into the de 
Sitter space $M^{n+1}_1$ and the projective space $P^{n+1}$, 
both of dimension $n+1$, and apply this mapping to find a 
set of their singular points. Theorems 4 and 5 are the core 
theorems of this paper. They give necessary and sufficient 
conditions for a submanifold  respectively in the de Sitter space 
$M^{n+1}_1$ and the projective space $P^{n+1}$ to be an image of 
an $m$-canal hypersurface of the conformal space $C^n$ 
under Darboux mapping. Moreover, Theorem 5 establishes the 
 very close connection between $m$-canal hypersurfaces 
in $C^n$ and tangentially degenerate submanifolds of $P^{n+1}$. 

Note that in this paper we 
consider  the real geometry of canal hypersurfaces.

{\bf 1.} {\bf The fundamental equations and notions of the 
conformal theory of hypersurfaces}. Let $V^{n-1}$ be a 
hypersurface of a conformal space $C^n$. With any point 
$x \in V^{n-1}$ we associate a family of conformal frames 
consisting of two points $A_0$ and $A_{n+1}$ and $n$ 
hyperspheres $A_i, i = 1, \ldots , n-1$, and $A_n$ passing through the points $A_0$ and $A_{n+1}$. We assume that the 
hypersphere $A_n$ is tangent to $V^{n-1}$ at the point $x$, and 
the hyperspheres $A_i$ are orthogonal to $V^{n-1}$ at $x$. 

Then the elements of the moving frame satisfy the following 
equations:
\begin{equation}\label{eq:1}
\left\{
\renewcommand{\arraystretch}{1.3}
\begin{array}{ll}
(A_0, A_0) = (A_{n+1}, A_{n+1}) = 0, \;\;
(A_{n+1}, A_i) = (A_{n+1}, A_n) = 0, \;\; (A_0, A_i) = 0, \\
(A_i, A_n) = 0, \;\; (A_i, A_j) = g_{ij}, \;\;
(A_n, A_n) = 1, \;\; (A_0, A_{n+1}) = -1, 
\end{array}
\renewcommand{\arraystretch}{1}
\right.
\end{equation}
where $(\;\; , \;\;)$ is the scalar product of the frame elements 
which is defined by the fundamental quadratic form of the space 
$C^n$:
\begin{equation}\label{eq:2}
(X, X) = g_{ij} x^i x^j + (x^n)^2 - 2 x^0 x^{n+1}.
\end{equation}
Here $x^0, x^i$, and $x^{n+1}$ are the polyspherical coordinates 
of elements (points and hyperspheres) of the space $C^n$. Note 
that $(X, X)= 0$ if $X$ is a point, $(X, X) > 0$ if $X$ is a 
hypersphere, and $(X, Y)= 0$ if $X$ is a point and $Y$ is a 
hypersphere passing through $X$ or $X$ and $Y$ are two orthogonal 
hyperspheres. The quadratic form $g_{ij} x^i x^j $ occurring in 
(2) is positive definite, and the signature of $(X, X)$ is 
$(n+1, 1)$. Note also that the last two relations in (1) are 
obtained as a result of  appropriate normalizations of the 
hypersphere $A_n$ and the points $A_0$ and $A_{n+1}$.

The equations of infinitesimal displacement of the conformal 
frame $\{A_\xi\}, \linebreak \xi = 0, 1, \ldots, n+1$, have the 
form
\begin{equation}\label{eq:3}
d A_\xi = \omega_\xi^\eta A_\eta, \;\;\;\; \xi, \eta 
= 0, 1, \ldots, n+1,
\end{equation}
where $\omega_\xi^\eta$ are 1-forms satisfying the structure 
equations of the space $C^n$:
\begin{equation}\label{eq:4}
d \omega_\xi^\eta = \omega_\xi^\zeta \wedge \omega_\zeta^\eta.
\end{equation}
Differentiating (1) by means of (3), we find that the forms 
$\omega_\xi^\eta$ satisfy the following Pfaffian equations:
\begin{equation}\label{eq:5}
\left\{
\renewcommand{\arraystretch}{1.3}
\begin{array}{ll}
\omega_0^{n+1} = 0, \;\; \omega^0_{n+1} = 0, \;\;  \omega_i^n 
+ g_{ij} \omega_n^j = 0, \;\; \omega_i^{n+1} = g_{ij} \omega_0^j,   
\;\; \omega_0^n - \omega_n^{n+1} = 0, \\
\omega_i^0 = g_{ij} \omega_{n+1}^j, \;\; 
\omega_n^0 - \omega_{n+1}^n = 0, \;\; dg_{ij} 
= g_{jk} \omega_i^k + g_{ik} \omega_j^k, \\
\omega_n^n = 0, \;\; \omega_0^0  + \omega_{n+1}^{n+1} = 0,
\end{array}
\renewcommand{\arraystretch}{1}
\right.
\end{equation}
where $\omega^i = \omega_0^i$ are basis forms of the hypersurface 
$V^{n-1}$. 

Since the hypersphere $A_n$ is tangent to $V^{n-1}$, we have 
$(dA_0, A_n) = 0$. By (3), this implies that 
\begin{equation}\label{eq:6}
\omega_0^n = 0.
\end{equation}
Equations (3) and (6) imply that 
\begin{equation}\label{eq:7}
dA_0 = \omega_0^0 A_0 + \omega^i A_i,
\end{equation}
where $\omega^i = \omega_0^i$ are basis forms of the hypersurface 
$V^{n-1}$. It follows from (7) that 
\begin{equation}\label{eq:8}
g = (dA_0, dA_0) = g_{ij} \omega^i \omega^j
\end{equation}
is a positive definite quadratic form determining a conformal 
structure on $V^{n-1}$. The quantities $g_{ij}$ are the 
components of a relative $(0, 2)$-tensor which is defined in a 
first-order differential neighborhood of a point $x \in V^{n-1}$.

It follows from (6) and (5) that 
\begin{equation}\label{eq:9}
\omega_n^{n+1} = 0.
\end{equation}
Taking the exterior derivative of equation (6), we obtain the 
exterior quadratic equation $\omega_i^n \wedge \omega^i = 0$. 
Applying Cartan's lemma to this equation, we find that 
\begin{equation}\label{eq:10}
\omega_i^n = \lambda_{ij} \omega^j, \;\; 
\lambda_{ij} = \lambda_{ji}.
\end{equation}
Taking exterior derivatives of equation (10) and applying 
Cartan's lemma to the exterior quadratic equation obtained, we 
arrive at the equations 
\begin{equation}\label{eq:11}
d \lambda_{ij} - \lambda_{kj} \omega_i^k 
- \lambda_{ik} \omega_j^k  + \lambda_{ij} \omega_0^0 
+ g_{ij} \omega_n^0 = \lambda_{ijk} \omega^k,    
\end{equation}
where $\lambda_{ijk}$ are symmetric with respect to all indices.

The quantities $\{g_{ij}\}, \{\lambda_{ij}, g_{ij}\}$, and 
$\{\lambda_{ijk}, \lambda_{ij}, g_{ij}\}$ form the geometric 
objects of $V^{n-1}$ which are connected with differential 
neighborhoods of a point $x \in V^{n-1}$ of orders 1, 2, and 3, 
respectively. 

The object $\{\lambda_{ij}, g_{ij}\}$ allows us to construct the 
following $(0, 2)$-tensor:
\begin{equation}\label{eq:12}
a_{ij} = \lambda_{ij} - \lambda g_{ij},
\end{equation}
where 
\begin{equation}\label{eq:13}
 \lambda = \frac{1}{n-1} g^{ij} \lambda_{ij}.
\end{equation}
The tensor $a_{ij}$ is apolar to the tensor $g_{ij}$: 
\begin{equation}\label{eq:14}
 a_{ij} g^{ij} = 0,
\end{equation}
where $g^{ij}$ is the inverse tensor of the tensor $g_{ij}$.

If $a_{ij}$ is a nondegenerate tensor (i.e., if $\det (a_{ij}) 
\neq 0$), then the object $\{\lambda_{ijk}, \lambda_{ij}, 
g_{ij}\}$ allows us to construct the following $(0, 3)$-tensor 
associated with a third-order differential neighborhood of a 
point $x \in V^{n-1}$: 
\begin{equation}\label{eq:15}
a_{ijk} = \lambda_{ijk} + 3(a_{(ij} g_{k)l} 
- g_{(ij} a_{k)l}) \mu^l,
\end{equation}
where 
\begin{equation}\label{eq:16}
\mu^i = a^{ij} \lambda_j, \;\;  \lambda_k = \frac{1}{n-1} g^{ij} \lambda_{ijk},
\end{equation}
and $a^{ij}$ is the inverse tensor of the tensor $a_{ij}$.

It is easy to check that the tensor $a_{ijk}$ is also apolar to 
the tensor $g_{ij}$: 
\begin{equation}\label{eq:17}
 a_{ijk} g^{ij} = 0. 
\end{equation}

We recall now the definition of the curvature lines on  a 
hypersurface  \linebreak $V^{n-1} \subset C^n$. Consider a 
symmetric affinor 
\begin{equation}\label{eq:18}
a^i_j = g^{ik} a_{kj},
\end{equation}
which is constructed by means of the tensors $g_{ij}$ and 
$a_{ij}$ and is called the {\em affinor of Burali-Forti} (see  
[Bu 12]). It is easy to see that the tensor $a^i_j$ is symmetric, 
and by (14), it is trace-free.

The directions that emanate from the point $A_0$ and are 
invariant with respect to the affinor $a_j^i$ (the 
eigendirections of  the affinor $a_j^i$) are called the 
{\em principal directions} of the hypersurface $V^{n-1}$. They 
are determined by the equations
\begin{equation}\label{eq:19}
a^i_j \xi^j = s \xi^i.
\end{equation}
The system (19) has a nontrivial solution if and only if 
\begin{equation}\label{eq:20}
\det (a_j^i - s \delta_j^i) = 0.
\end{equation}
Equation (20) is equivalent to the equation 
\begin{equation}\label{eq:21}
\det (a_{ij} - s g_{ij}) = 0.
\end{equation}
Equation (20) is the characteristic equation of the affinor 
$a_j^i$. Its roots are the eigenvalues of the affinor $a_j^i$. 
Since the affinor $a_j^i$ is symmetric, its eigenvalues are real, 
and, by (18) and (14),  their sum is equal to 0. 

If all eigenvalues of the affinor $a_j^i$ are mutually distinct, 
then the corresponding principal directions are orthogonal. The 
curves of $V^{n-1}$ enveloping the principal directions are 
called the {\em curvature lines}. If this is the case, we say 
that {\em the hypersurface $V^{n-1}$ carries a unique net of 
curvature lines.} 

If the affinor $a_j^i$ has an eigenvalue of multiplicity $m > 1$, 
then it defines an $m$-dimensional domain $L^m$ of principal 
directions, and each direction of $L^m$ is a principal direction.

Consider a hypersphere $A_n + s A_0$ that is tangent to the 
hypersurface $V^{n-1}$ at its point $A_0$. Let us find those 
directions on  $V^{n-1}$ along which this hypersphere has a 
second-order tangency with  $V^{n-1}$. Such directions are 
determined by the equation
$$
(d^2 A_0, A_n + s A_0) 
= (\lambda_{ij} - s g_{ij}) \omega^i \omega^j = 0.
$$
Thus the desired  directions constitute a cone of second order 
with its vertex at the point $A_0$. In general, the rank of the 
quadratic form $(a_{ij} - s g_{ij}) \omega^i \omega^j$ is equal 
to $n-1$, since the  tensor $g_{ij}$ is nondegenerate. However, 
for some values  of $s$, this rank can be reduced. These values 
of $s$ are determined by the equation 
\begin{equation}\label{eq:22}
\det (\lambda_{ij} - s g_{ij}) = 0.
\end{equation}
Essentially equation (22) coincides with equation (20). Suppose 
that $s_h$ are the roots of equation (22). Then the hypersphere 
$B_h = A_n + s_h A_0$ has a second-order tangency with the 
hypersurface $V^{n-1}$ along the directions defined by the 
equation
\begin{equation}\label{eq:23}
(\lambda_{ij} - s_h g_{ij}) \omega^i \omega^j = 0. 
\end{equation}
The rank of the quadratic form in the left-hand side of 
equation (23) is equal to $n - m - 1$, where $m$ is the 
multiplicity of the root $s_h$ of the characteristic equation 
(22). Thus  equation (23) defines a degenerate cone with an 
$m$-dimensional vertex. This vertex is defined by the system of 
equations
\begin{equation}\label{eq:24}
(\lambda_{ij} - s_h g_{ij}) \omega^j = 0 
\end{equation}
and coincides with the $m$-dimensional eigendirection $L^m$ 
corresponding to the $m$-multiple root $s_h$. If $s_h$ is a 
simple root of equation (22), then the vertex of the 
corresponding cone is one-dimensional.

{\bf 2.} {\bf Analytic characterization of  
\protect\boldmath\(m\)\protect\unboldmath-canal hypersurfaces for 
\protect\boldmath\(m \geq \;\)\protect\unboldmath \nolinebreak 
2}. We now 
consider a hypersurface $V^{n-1}$ whose affinor $a^j_i$ has an  
eigenvalue of multiplicity $m$ at each point $A_0 \in V^{n-1}$. 
Suppose, for example, that all eigenvalues of the affinor $a^j_i$ 
satisfy the conditions
\begin{equation}\label{eq:25}
\renewcommand{\arraystretch}{1.3}
\begin{array}{ll}
a_1 = a_2 = \ldots = a_m = a, \;\; m \geq 2, \;\;
a\neq a_p,  \;\;\;\;\; p = m+1, \ldots , n - 1.
\end{array}
\renewcommand{\arraystretch}{1}
\end{equation}
 Then the rank of the tensor $a_i^j - a \delta_i^j$ is equal to 
$r = n - m - 1$, and at each point $A_0$ the system of equations 
(19) for $s = a$ determines an $m$-dimensional subspace $L^m$ of 
principal directions corresponding to this $m$-multiple root $a$. 

If a hypersurface $V^{n-1}$ is the envelope of an $r$-parameter 
family of hyperspheres and $m = n - r -1$, it  is called an 
{\em $m$-canal hypersurface}. Such a hypersurface carries an 
$r$-parameter family of $m$-dimensional spherical characteristics 
$S^m$. 

We will now prove   the following theorem:

\begin{theorem}
For $m \geq 2$, a hypersurface $V^{n-1}$ is an $m$-canal 
hypersurface if and only if its tensor  $a_{ij}$ has an 
eigenvalue of multiplicity  $m$. 
\end{theorem}

{\sf Proof.} First we will prove the {\em sufficiency}. Let $L^m$ 
be a subspace of principal directions corresponding to the root 
$a$ of multiplicity $m$ of the affinor $a^i_j$ of the 
hypersurface $V^{n-1}$. With a point $A_0 \in V^{n-1}$ we 
associate a frame bundle in such a way that its hyperspheres 
$A_a, \; a = 1, \ldots , m$,  are orthogonal to the subspace 
$L^m$, and the hyperspheres $A_p, \; p = m+1, \ldots , n-1$, are 
tangent to $L^m$. Then, in the chosen frame bundle we have
\begin{equation}\label{eq:26}
(g_{ij}) = \pmatrix{g_{ab} & 0 \cr
                  0      & g_{pq}}, \;\;  
(\lambda_{ij}) = \pmatrix{\widetilde{\lambda} g_{ab} & 0 \cr
                  0      & \lambda_{pq}}, \;\;  
(a_{ij}) = \pmatrix{a g_{ab} & 0 \cr
                  0      & a_{pq}}, 
\end{equation}
where $\widetilde{\lambda} = a + \lambda$ and $a = - \frac{1}{m} 
a_{ab} g^{ab}$.

The hypersphere $B = A_n + \widetilde{\lambda} A_0$ is a contact 
hypersphere for the subspace $L^m$. In fact, we have 
$$
(d^2 A_0, B) = (\lambda_{ij} - \widetilde{\lambda} g_{ij}) 
\omega^i \omega^j = (\lambda_{pq} - \widetilde{\lambda} g_{pq}) 
\omega^p \omega^q,
$$
and as a result, the hypersphere $B$ has a second order 
tangency along the direction $L^m$ defined by the equations 
$\omega^p = 0$. 

We specialize our frame bundle by taking $A_n = B$. Then we have 
$\widetilde{\lambda} = 0$ and 
\begin{equation}\label{eq:27}
(\lambda_{ij}) = \pmatrix{0 & 0 \cr
                  0      & \lambda_{pq}},
\end{equation}
where $(\lambda_{pq})$ is a nondegenerate matrix

We will now write the differential equations which the quantities 
$g_{ij}$ and $\lambda_{ij}$ satisfy in our new specialized 
frame bundle. It follows from (5) that 
\begin{equation}\label{eq:28}
\left\{
\renewcommand{\arraystretch}{1.3}
\begin{array}{ll}
dg_{ab} = g_{ac} \omega_b^c + g_{cb} \omega_a^c, \\
g_{ab} \omega_p^b + g_{pq} \omega_a^q = 0, \\
dg_{pq} = g_{pt} \omega_q^t + g_{tq} \omega_a^p,
\end{array}
\renewcommand{\arraystretch}{1}
\right.
\end{equation}
and it follows from (11) that 
\begin{equation}\label{eq:29}
\left\{
\renewcommand{\arraystretch}{1.3}
\begin{array}{ll}
g_{ab} \omega_n^0 = \lambda_{abk} \omega^k, \;\;
\lambda_{pq} \omega_a^q = - \lambda_{apk} \omega^k, \\
d \lambda_{pq} - \lambda_{pt} \omega_q^t - \lambda_{tq}\omega_p^t  + \lambda_{pq} \omega_0^0 + g_{pq} \omega_n^0 = \lambda_{pqk} 
 \omega^k,
\end{array}
\renewcommand{\arraystretch}{1}
\right.
\end{equation}
where here and in what follows the index $t$ has the same range 
as the indices $p$ and $q$: $t = m+1, \ldots , n-1$. 

It follows from the first equation of (29) that 
\begin{equation}\label{eq:30}
 \omega_n^0 = b_k \omega^k. 
\end{equation}
Substituting this value into the first equation of (29), we 
find that 
\begin{equation}\label{eq:31}
 \lambda_{abk} = g_{ab} b_k. 
\end{equation}
But the quantities $\lambda_{ijk}$ are symmetric. Thus 
$$
\lambda_{abc} = g_{ab} b_c = g_{bc} b_a = g_{ca} b_b.
$$
Contracting the middle equation with the tensor $g^{ab}$, we 
obtain
$$
m b_c = b_c, \;\;\;\;\; c = 1, \ldots , m.
$$
Since we assume that $m \geq 2$, it follows that 
$$
b_c = 0,  \;\;\;\;\; c = 1, \ldots , m. 
$$
Thus it follows from equation (31) that 
\begin{equation}\label{eq:32}
\lambda_{abc} = 0, \;\;  \lambda_{abp} = g_{ab} b_p, 
\end{equation}
and equation (30) takes the form
\begin{equation}\label{eq:33}
 \omega_n^0 = b_p \omega^p. 
\end{equation}
Note that the above computation is the key step in the proof of 
Theorem 1. For $m = 1$ a relation similar to the first equation of (32) is not valid. 

Now the second equation of system (29) can be written as 
\begin{equation}\label{eq:34}
d \lambda_{pq} - \lambda_{pt} \omega_q^t - \lambda_{tq}\omega_p^t  + \lambda_{pq} \omega_0^0 = \lambda_{pqa}  \omega^a 
+ (\lambda_{pqt} - g_{pq} b_t) \omega^t. 
\end{equation}
This shows that the quantities $\lambda_{pq}$ form a tensor with 
respect to admissible transformations of specialized frames. By 
(27), this tensor is nondegenerate since the rank of the matrix 
$(\lambda_{ij})$ is equal to $r$. Hence we can write the second 
equation of (29) in the form 
\begin{equation}\label{eq:35}
 \omega_a^p = -  \lambda^{pq} (g_{ab} b_q \omega^b 
+ \lambda_{aqt} \omega^t).
\end{equation}

Next we consider the family of contact hyperspheres $A_n$. 
Differentiating $A_n$ and applying equations (5), we find that 
$$
dA_n = \omega_n^0 A_0 + \omega_n^a A_a + \omega_n^p A_p.
$$
But the form $\omega_n^0$ is expressed by formula (33), 
and by (27) we have 
\begin{equation}\label{eq:36}
 \omega_n^a = 0, \;\; \omega_n^p = -\lambda_q^p \omega^q,
\end{equation}
where $\lambda^p_q = g^{pt} \lambda_{tq}$ is a nondegenerate 
affinor. Thus we have 
\begin{equation}\label{eq:37}
dA_n = (b_p A_0 - \lambda_p^q A_q) \omega^p.
\end{equation}

Hence the hypersphere $A_n$ depends on $r = n-m-1$ parameters. 
Since by (35) we have
$$
d \omega^p \equiv 0 \pmod{\omega^p},
$$
the system of equations
$$
\omega^p = 0
$$
is completely integrable on  the hypersurface $V^{n-1}$. On 
$V^{n-1}$ these equations define $m$-dimensional characteristics 
of the family of hyperspheres $A_n$, and along these 
characteristics the hyperspheres $A_n$ are tangent to the 
hypersurface $V^{n-1}$. These characteristics  are the 
intersections of the hypersphere $A_n$ with the linearly 
independent hyperspheres 
\begin{equation}\label{eq:38}
C_p = b_p A_0 - \lambda_p^q A_q.
\end{equation}
and they are $m$-dimensional spherical generators of the 
hypersurface $V^{n-1}$. 

Thus, if $m$ eigenvalues of the tensor $a_{ij}$ of the 
hypersurface $V^{n-1}$ coincide, where $m \geq 2$, then this 
hypersurface is  an  $m$-canal hypersurface.
 
Let us prove the {\em necessity} of the theorem conditions:  
if a hypersurface $V^{n-1}$ is an  $m$-canal hypersurface where 
$m \geq 2$, then its tensor $a_{ij}$ has an eigenvalue of 
multiplicity $m$. In fact, an  $m$-canal hypersurface is the 
envelope of an $r$-parameter family of contact  hyperspheres $B$, 
where $r = n - m - 1$. By specializing our moving frame, we make 
its hypersphere $A_n$ coincide with the hypersphere $B$. Since 
\begin{equation}\label{eq:39}
dA_n = \omega_n^0 A_0 + \omega_n^i A_i,
\end{equation}
the forms $\omega_n^0$ and $\omega_n^i$ must be expressed in 
terms of $r$ linearly independent forms $\omega^p$. But we have 
$$
\omega_n^i = - g^{ij} \lambda_{jk} \omega^k.
$$
Thus the rank of the matrix $(\lambda_{ij})$ is equal to $r$, and 
by (12), the tensor $a_{ij}$ has precisely $m$ equal eigenvalues. 
This concludes the proof of Theorem 1.

Note that if  $r = 0$, then all eigenvalues of the tensor  
$a_{ij}$ coincide, and by  the apolarity condition (14),  all 
these eigenvalues are equal to 0. It follows that $a_{ij} = 0$ 
and $\omega_n^i = 0$. Exterior differentiation of the last 
equations leads to the exterior quadratic equations 
$$
\omega_n^0 \wedge \omega_0^i = 0,
$$
where $i = 1, \ldots , n - 1$ and $n - 1 \geq 2$. 
Thus the last exterior equations imply that $\omega_n^0 = 0$. By 
(39), this means that the  hypersphere is fixed, and the 
hypersurface $V^{n-1}$ coincides with the  hypersphere $A_n$ or 
with its open part.

 {\bf 3.}  {\bf Analytic characterization of 1-canal 
hypersurfaces}. As we saw in Section {\bf 2}, for $m \geq 2$, 
the analytic characterization of $m$-canal hypersurfaces is 
expressed in terms of the quantities connected with a 
second-order neighborhood of a point $A_0 \in V^{n-1}$. In 
contrast, for $m =1$, the analytic characterization of $1$-canal 
hypersurfaces is expressed in terms of the quantities connected 
with a third-order neighborhood of a point $A_0 \in V^{n-1}$.

Namely the following theorem is valid:

\begin{theorem}
Suppose that $s_1$ is a simple eigenvalue of the affinor $a_j^i$ 
of the 
hypersurface $V^{n-1}$ and  that $V^{n-1}$ 
is referred to a frame whose hypersphere $A_1$ is 
orthogonal to the corresponding eigendirection 
of this affinor and whose hyperspheres $A_p, p \neq 1$ 
are tangent to this direction.  Then the hypersurface 
$V^{n-1}$ is a $1$-canal hypersurface for this 
eigendirection if and only if 
in the above frame the diagonal component $a_{111}$ of the tensor 
$a_{ijk}$ vanishes. 
\end{theorem}

{\sf Proof}. On the hypersurface $V^{n-1}$ we consider a family 
of contact hyperspheres $B_1 = A_n + s_1 A_0$, where $s_1$ is a 
simple 
root of characteristic equation (22). We specialize a family of 
frames associated with a point $A_0 \in V^{n-1}$ in such a way 
that the tangent hypersphere $A_n$ coincides with the contact 
hypersphere $B_1$. Then $s_1 = 0$, and the hypersphere $A_n$ has 
a second-order tangency with the hypersurface $V^{n-1}$ along the 
direction determined by the equations $\omega^p = 0, \; p = 2, 
\ldots , n - 1$. If in addition we specialize  our 
moving frame as indicated in the theorem condition, then the 
matrices $(g_{ij}), (\lambda_{ij})$, and $(a_{ij})$ take the form 
\begin{equation}\label{eq:40}
(g_{ij}) = \pmatrix{1 & 0 \cr
                  0      & g_{pq}}, \;\;  
(\lambda_{ij}) = \pmatrix{0 & 0 \cr
                  0      & \lambda_{pq}}, \;\;  
(a_{ij}) = \pmatrix{a  & 0 \cr
                  0    & a_{pq}}, 
\end{equation}
where $p, q = 2, \ldots , n - 1; \; g_{pq}, \lambda_{pq}$, and 
$a_{pq}$ are nondegenerate tensors, and \linebreak 
$a = - g^{pq} a_{pq}$. 

Now those of equations (5), which the components of the tensor 
$g_{ij}$ satisfy, and also equations (11) take the form
\begin{equation}\label{eq:41}
\left\{
\renewcommand{\arraystretch}{1.3}
\begin{array}{ll}
\omega_1^1 = 0, \;\; \omega^1_p + g_{pq} \omega_1^q = 0, \\  
dg_{pq} = g_{pt} \omega_q^t + g_{tq} \omega_p^t, 
\end{array}
\renewcommand{\arraystretch}{1}
\right.
\end{equation}
and 
\begin{equation}\label{eq:42}
\left\{
\renewcommand{\arraystretch}{1.3}
\begin{array}{ll}
\omega_n^0 = \lambda_{11k} \omega^k, \;\; 
\lambda_{pq} \omega_1^q = - \lambda_{1pk} \omega^k, \\
d \lambda_{pq} - \lambda_{pt} \omega_q^t 
- \lambda_{tq} \omega_p^t  + \lambda_{pq} \omega_0^0 
+ g_{pq} \omega_n^0 = \lambda_{pqk} \omega^k,    
\end{array}
\renewcommand{\arraystretch}{1}
\right.
\end{equation}
where $p, q, t = 2, \ldots , n-1$.

Next we find the differential of the contact hypersphere 
$A_n$ when it is moving along the hypersurface $V^{n-1}$:
$$
dA_n = \omega_n^0 A_0 + \omega_n^1 A_1 + \omega_n^p A_p.
$$
The form $\omega_n^0$ occurring in the last equation is expressed 
by formula (42), and by (40), we have 
$$
\omega_n^1 = 0, \;\; \omega_n^p = - \lambda^p_q \omega^q,
$$
where $\lambda_q^p = g^{pt} \lambda_{tq}$ is a 
nondegenerate affinor. This implies that 
$$
dA_n = \lambda_{111} \omega^1 A_0 + 
(\lambda_{11p} A_0 - \lambda^q_p A_q)\omega^p.
$$
By (40), it follows from (15) that $\lambda_{111} = a_{111}$. 
Thus the last equation can be written in the form
\begin{equation}\label{eq:43}
dA_n = a_{111} \omega^1 A_0 + 
(\lambda_{11p} A_0 - \lambda^q_p A_q)\omega^p.
\end{equation}
This shows that the hypersphere $A_n$ depends on $n-2$ 
parameters if and only if $a_{111} = 0$. In this case the 
hypersurface $V^{n-1}$ is the envelope of an $(n-2)$-parameter 
family of hyperspheres $A_n$, that is, the hypersurface $V^{n-1}$ 
is a 1-canal hypersurface. This completes the proof of Theorem 2. 
 
Characteristics of the $(n-2)$-parameter family of contact 
hyperspheres $A_n$ are circles $S^1$ which are generators of 
the  hypersurface $V^{n-1}$.  These circles form the first 
family of curvature lines of the   hypersurface $V^{n-1}$. 

Note that relation (43) differs from relation (37) 
which is valid for $m \geq 2$ by the presence of the first term 
containing the coefficient $a_{111}$. Namely this distinguishes 
the cases $m = 1$ and $m \geq 2$. For a 1-canal hypersurface we 
have $a_{111} = 0$. If in addition to this condition we set 
$\lambda_{11p} = b_p$ in equation (43), then  equation (43) 
completely coincides with  equation (37). We will use equation 
(37) for any $m = 1, 2, \ldots , n - 2$. This allows us 
to consider the cases $m = 1$ and $m \geq 2$ simultaneously. 

Theorem 2 implies the following method for determining whether a 
hypersurface $V^{n-1}$, whose affinor $a_j^i$ has a simple 
eigenvalue $s_i$, is a 1-canal hypersurface. To this end, 
we must choose a frame subbundle associated with $V^{n-1}$ 
in such a way that its hypersphere $A_i$ is orthogonal to 
the corresponding eigendirection and compute the tensor 
$a_{ijk}$ in this frame. If its diagonal component 
$a_{iii}$ vanishes on $V^{n-1}$, then the hypersurface 
$V^{n-1}$ is a 1-canal hypersurface for the direction 
$\omega^k = 0, k \neq i$.

If on  a hypersurface $V^{n-1}$ the conditions $a_{iii} = 0$ hold 
for all values of $i$, then this hypersurface is the envelope of 
$n-1$  families of hyperspheres, each of which depends on $n-2$ 
parameters. Such  a hypersurface is a multidimensional analog of 
the well-known Dupin cyclide (see [P 85], [V 58], and the book 
[C 92] where one can find a detailed bibliography on Dupin's 
cyclide).

Since the tensor $a_{ijk}$ is apolar to the tensor $g_{ij}$, in 
our moving frame we find that if $n = 3$, then the condition 
$a_{iii} = 0$ implies that $a_{iij} = 0$; that is, the conditions 
$a_{ijk} = 0$ hold for any values of the indices $i, j$, and $k$. 
Since $a_{ijk}$ is a tensor, it vanishes not only in the special 
frame under consideration but also in any first-order frame. 

 Thus we have proved the following theorem:

\begin{theorem}
A two-dimensional surface $V^2$ of a three-dimensional conformal 
space $C^3$ is a Dupin cyclide if and only if its tensor 
$a_{ijk}$, determined by a third-order differential neighborhood, 
vanishes.
\end{theorem}

This differential geometric characterization of the Dupin cyclide 
is of invariant nature; that is, it does not depend on either 
the choice of a coordinate system on the surface $V^2$ or the 
choice of a conformal frame associated with the surface $V^2$.
Since the  tensor  $a_{ijk}$ can also be calculated for a surface 
$V^2$  of a three-dimensional Euclidean space, {\em the invariant  
characterization of the Dupin cyclides we have obtained is also 
valid for the  Dupin cyclides in the Euclidean space $R^3$}.  

{\bf 4. Submanifolds in the de Sitter space}. 
In Sections {\bf 4} and {\bf 5} we will not distinguish the cases 
$m = 1$ and $m \geq 2$. Consider the Darboux mapping of the 
conformal space $C^n$ (see [K 72] or [N 50], \S\S 28-29 or 
[AG 96], Sections {\bf 1.1} and {\bf 1.3}). Under this mapping to 
a point of $C^n$ there corresponds a point of an oval 
hyperquadric $Q^n$ of the projective space $P^{n+1}$. The 
hyperquadric $Q^n$ is called the {\em Darboux hyperquadric}. To a 
hypersphere of  $C^n$ there corresponds a point of $P^{n+1}$ 
lying outside of the Darboux hyperquadric $Q^n \subset P^{n+1}$. 
To the group of conformal transformations of the space $C^n$ 
there corresponds the group of projective transformations of the  
space $P^{n+1}$ whose transformations map the hyperquadric 
$Q^n$ onto itself. To any geometric property of the preimages of 
the space $C^n$ there corresponds a certain property of the 
images  of the space $P^{n+1}$ connected with the hyperquadric 
$Q^n$. 

A domain  outside of the Darboux hyperquadric $Q^n \subset 
P^{n+1}$, onto which the hyperspheres of the  space $C^n$ are 
mapped, 
is a pseudo-Riemannian manifold of index 1 and of constant 
positive curvature. This manifold is called the 
{\em de Sitter space} 
and is denoted by $M_1^{n+1}$ (see, for example, 
[Z 96]). The hyperquadric $Q^n$ is the absolute of the space 
$M_1^{n+1}$.

Let $x$ be an arbitrary point of the space $M_1^{n+1}$. The 
tangent lines to the hyperquadric $Q^n$ passing through the 
point $x$ form a cone $C_x$ with vertex at $x$. This cone 
is called the {\em isotropic cone}. The cone $C_x$ separates all 
straight lines passing through the point $x$ into spacelike 
straight lines that do not have points in common with the 
hyperquadric $Q^n$, timelike straight lines intersecting $Q^n$ at 
two different points, and lightlike straight lines that are 
tangent to the hyperquadric $Q^n$. The lightlike straight lines 
are generators of the cone $C_x$. This is the reason that the 
cone $C_x$ is also called the {\em light cone}.

To a spacelike straight line $l \subset M_1^{n+1}$ there 
corresponds an elliptic pencil of hyperspheres in the conformal 
space $C^n$. All hyperspheres of this pencil pass through a 
common $(n-2)$-sphere $S^{n-2}$ (the center of this pencil). 
The sphere $S^{n-2}$ is the intersection of the hyperquadric 
$Q^n$ and an $(n-1)$-dimensional subspace of the space $P^{n+1}$ 
which is polar conjugate to the line $l$ with respect to the 
hyperquadric $Q^n$. To a timelike straight line $l \subset 
M_1^{n+1}$ there corresponds a hyperbolic pencil of hyperspheres 
in the space $C^n$. Two arbitrary hyperspheres of this pencil do 
not have common points, and the pencil contains two hyperspheres 
of zero radius which correspond to the points of intersection of 
the straight line $l$ and the hyperquadric $Q^n$. Finally, to a 
lightlike straight line $l 
\subset M_1^{n+1}$ there corresponds a parabolic pencil of 
hyperspheres in the  space $C^n$ consisting of 
hyperspheres tangent one to another at a point that is a unique 
hypersphere of zero radius belonging to this pencil.

The curves in the space $M_1^{n+1}$ are also divided into three 
classes: spacelike, timelike, and lightlike. It is not difficult 
to prove (see [AG 97]) that lightlike curves are always straight 
lines, and that to these lines there corresponds a parabolic 
pencil of hyperspheres in $C^n$. To a spacelike curve $x = x (t)$ 
of the space  $M_1^{n+1}$, for which $x' (t) \neq 0$ and $x'' (t) \neq 0$, 
 there corresponds a one-parameter 
family of hyperspheres in $C^n$ having an $(n-2)$-canal 
hypersurface as its envelope. Characteristics of this canal 
hypersurface are $(n-2)$-dimensional spheres that correspond to 
 spacelike tangents to the curve $x (t)$. Finally, to a 
lightlike curve $x = x (t)$ of the space  $M_1^{n+1}$, 
there corresponds a one-parameter family of hyperspheres in $C^n$ not 
having an envelope. In fact, neighboring hyperspheres of this 
family that are determined by timelike tangents to the curve 
$x (t)$ do not have common points, and as a result, the family 
does not possess characteristics. 

In the de Sitter space  $M_1^{n+1}$ we consider a submanifold 
$U^r$ of dimension $r, \;1 \leq r \leq n - 2$. Such submanifolds 
are also 
divided into three classes:  spacelike, timelike, and lightlike. 
The tangent subspaces to spacelike submanifolds do not have 
common points with the Darboux hyperquadric $Q^n$, the tangent 
subspaces to timelike submanifolds intersect $Q^n$ along an 
$(r-1)$-dimensional sphere, and the tangent subspaces to 
lightlike submanifolds are tangent to $Q^n$.

Consider now an $m$-canal hypersurface in the conformal space 
$C^n$. By definition, such a hypersurface is the envelope of an 
$r$-parameter differentiable family of hyperspheres $S^{n-1}$, 
where $r = n - m - 1$. To such a family there corresponds an 
$r$-dimensional differentiable submanifold $U^r$ in the space 
$M^{n+1}_1$. However, this submanifold is not arbitrary. Namely, 
the following theorem is valid:

\begin{theorem} A differentiable submanifold $U^r$ of dimension 
$r$ in the space $M^{n+1}_1$ corresponds to an $m$-canal 
hypersurface $V^{n-1}, \; m = n-r-1$, in the conformal space 
$C^n$ if and only if $U^r$ is tangentially nondegenerate and 
spacelike. 
\end{theorem}

{\sf Proof}. In fact, let  $V^{n-1}$ be an $m$-canal hypersurface 
in $C^n, m = 1, \ldots , n-2$. As we showed in Sections {\bf 2} 
and {\bf 3}, for a set of contact hyperspheres of $V^{n-1}$ whose 
family is enveloped by the hypersurface $V^{n-1}$ itself,  
condition (37) is satisfied. This condition can be written in the 
form
$$
dA_n = \omega^p C_p,
$$
where the hyperspheres $C_p$ are defined by formula (38). 
In the de Sitter space $M^{n+1}_1$, to an $m$-canal hypersurface 
$V^{n-1}$ there corresponds a submanifold $U^r$ described 
by the point $A_n$ corresponding to the hypersphere $A_n 
\subset C^n$ whose tangent subspace is determined by 
the points $A_n$ and $C_p$. An arbitrary point $z$ of 
this tangent subspace can be written in the form
$$
z =z^n A_n + z^p C_p = z^n A_n + z^p (b_p A_0 - \lambda_p^q A_q). 
$$
Substituting the coordinates of this point into equation (2), 
we find that 
$$
(z, z) = g_{st} \lambda_i^s \lambda_q^t z^p z^q + (z^n)^2 > 0. 
$$
It follows that the tangent subspace to the submanifold $U^r$ 
at the point $A_n$ have no common points with the 
hyperquadric $Q^n$, that is, the tangent subspace to $U^r$ 
at $A_n$ is spacelike. 
Since this is true for any point $A_n \in U^r$,   the submanifold 
$U^r$ corresponding to the family of hyperspheres $A_n$ in the 
space $M^{n+1}_1$ is spacelike. 

The tangential nondegeneracy  of 
the  submanifold $U^r$ follows from the fact that the tangent 
subspaces $T_{A_n} (U^r)$ of such a submanifold correspond to the 
characteristic spheres $S^m$ of the hypersurface $V^{n-1} 
\subset C^n$ that depend on $r$ parameters, and as a result, the 
tangent subspaces $T_{A_n} (U^r)$ depend on the same number of 
parameters.

If a submanifold $U^r$ in the space $M^{n+1}_1$ is timelike or 
lightlike, then it contains curves which correspond to 
one-parameter families of hyperspheres in $C^n$ not possessing 
envelopes. Thus an $r$-parameter family of hyperspheres in $C^n$, 
whose image in $M^{n+1}_1$ is the submanifold $U^r$, also does not 
possess an envelope, and hence such a family does not define a 
canal hypersurface.

{\bf 5. The Darboux mapping of 
\protect\boldmath\(m\)\protect\unboldmath-canal hypersurfaces 
\protect\boldmath\(V^{n-1} \subset C^n\)\protect\unboldmath.} 
In this section we will consider the same Darboux mapping which 
was considered in Section {\bf 4.} However, since here we 
will be interested in projective properties of this mapping 
and  its images, we will consider the mapping of the 
conformal space $C^n$ onto a projective space $P^{n+1}$ instead 
of the mapping of $C^n$ onto the de Sitter space $M^{n+1}_1$.

A canal hypersurface $V^{n-1}$ is the envelope of an 
$r$-parameter family of contact hyperspheres $B$, where 
$r = n - m - 1, r \leq n-2$. As we showed in Sections {\bf 2} 
and {\bf 3}, by specializing moving frames we can make the 
hypersphere $A_n$ coincide with the hypersphere $B$. 

Under the Darboux mapping, to the hypersphere $A_n$ there 
corresponds a point $x$ of the projective space $P^{n+1}$ lying 
outside of the Darboux hyperquadric $Q^n \subset P^{n+1}$ onto 
which the points of the conformal space $C^n$ are mapped 
bijectively. The polar hyperplane of the point $x$ with respect 
to  $Q^n$ is a hyperplane $\xi$ intersecting $Q^n$ at  
points which are images of points of the hypersphere $A_n \subset 
C^n$. The point $x$ describes in $P^{n+1}$ a tangentially 
nondegenerate spacelike submanifold $U^r$, and its polar 
hyperplane $\xi$ 
envelopes a tangentially degenerate hypersurface $U^n$ of rank 
$r$ (see [A 57] or [AG 93], Ch. 4). The  hypersurface $U^n$ 
carries timelike plane generators $\alpha^{m+1}$ of dimension 
$m+1$ which are polar conjugate 
to the $r$-dimensional tangent subspaces $T_x (U^r)$ of 
the submanifold $U^r$. At each point of a generator 
$\alpha^{m+1}$ the tangent hyperplane to the  hypersurface $U^n$ 
coincides with the hyperplane $\xi$. Thus we prove the following 
result:

\begin{theorem} 
The hypersurface $U^n$ of the projective space $P^{n+1}$, that is 
the image of an $m$-canal hypersurface $V^{n-1}$ under the 
Darboux mapping, is a tangentially degenerate hypersurface 
of rank $r = n - m - 1$ with $(m+1)$-dimensional plane 
generators. 
\end{theorem}
 
A plane generator $\alpha^{m+1}$ of a tangentially degenerate 
hypersurface $U^n$ carries a focus surface ${\cal F}, \dim {\cal 
F} = m$, formed by singular points of this generator (see [A 57] 
and [A 87]). The submanifold ${\cal F}$ is an algebraic 
hypersurface of degree $r$ in the generator $\alpha^{m+1}$ of the 
hypersurface $U^n$ that corresponds to a canal hypersurface 
$V^{n-1} \subset C^n$ under the Darboux mapping.

We will find the equation of the focus surface ${\cal F}$ of the 
generator $\alpha^{m+1}$ of a hypersurface $U^n \subset P^{n+1}$ 
that corresponds to an $m$-canal hypersurface of the conformal 
space $C^n$ under the Darboux mapping. As we already noted above, 
for canal hypersurfaces $V^{n-1}$ the differential (37) and (43) 
of contact hyperspheres will have the same form if we set 
$\lambda_{11p} = b_p$ in equations (43). This is the reason that 
we can consider these two cases simultaneously. 

Relations (37) are preserved under the Darboux mapping. Since the 
point $x \in U^r$ coincides with the point $A_n$, then by (37), 
the tangent subspace $T_x (U^r)$ to this submanifold $U^r$ is 
determined by the points $A_n$ and $C_p = b_p A_0 - \lambda_p^q A_q$ of the space $P^{n+1}$. If $x = x^0 A_0 + x^a A_a + x^p A_p 
+ x^n A_n + x^{n+1} A_{n+1}$ is an arbitrary point of the plane 
generator $\alpha^{m+1}  \subset U^n$ that is polar conjugate to 
the subspace $T_x (U^r)$, then the generator $\alpha^{m+1}$ is 
defined by the equations 
$$
(x, A_n) = x^n = 0, \;\; (x, C_p) 
= - b_p x^{n+1} - \lambda_{pq} x^q = 0.
$$
Thus an arbitrary point $x$ of the $\alpha^{m+1}$ can be 
expressed as follows:
$$
x = x^0 A_0 + x^a A_a + x^{n+1} (A_{n+1} - \lambda^{pq} b_q A_p), 
$$
where $\lambda^{pq}$ is the inverse tensor of the tensor 
$\lambda_{pq}$.

We can further specialize the frame bundle associated with 
the submanifold $U^n \subset P^{n+1}$ if we make the point $A_{n+1}$ coincide with the point $A_{n+1} - \lambda^{pq} b_q 
A_p$ located in the generator $\alpha^{m+1}$. Then we obtain
$$
b_p = 0,
$$
and the points $A_0, A_a$ and $A_{n+1}$ become the basis points 
of the generator $\alpha^{m+1}$. Now it follows from equation 
(32)  that 
$$
\lambda_{abk} = 0,
$$
and as a result, the first equation of (29) takes the form 
\begin{equation}\label{eq:44}
 \omega_n^0 =0.
\end{equation}

Now the second equation of (29) gives 
$$
\lambda_{pq} \omega_a^q = - \lambda_{apq} \omega^q, 
$$
and this implies that 
\begin{equation}\label{eq:45}
 \omega_a^p = - \lambda^{pt} \lambda_{atq} \omega^q 
= \lambda^p_{aq} \omega^q, 
\end{equation}
where $\lambda^p_{aq} = - \lambda^{pt} \lambda_{atq}$.

Exterior differentiation of equation (44) and application of (27) 
gives
$$
\lambda_p^q \omega_q^0 \wedge \omega^p_0 = 0.
$$
The last equation implies that 
\begin{equation}\label{eq:46}
 \omega_p^0 = c_{pq}  \omega^q,
\end{equation}
and the coefficients $c_{pq}$ satisfy the conditions:
\begin{equation}\label{eq:47}
\lambda_p^t c_{tq} = \lambda_q^t c_{tp}, 
\end{equation}
which can be obtained if one substitutes decompositions (46) into 
preceding exterior equations. 

Differentiating the basis points  the points $A_0, A_p$ and   
$A_{n+1}$ of the generator $\alpha^{m+1}$  and applying equations 
(5), (6),  and (44), we find that 
\begin{equation}\label{eq:48}
\renewcommand{\arraystretch}{1.3}
\begin{array}{lll}
dA_0 \!\!\!\!&= \omega_0^0 A_0 +& \omega_0^a A_a + \omega_0^p A_p, \\
dA_a \!\!\!\!&= \omega_a^0 A_0 +& \omega_a^b A_b + \omega_a^p A_p
                              +  \omega_a^{n+1} A_{n+1}, \\
dA_{n+1}\!\!\!\! &=         & \omega_{n+1}^a A_a + \omega_{n+1}^p A_p 
                                 - \omega_0^0 A_{n+1}.
\end{array}
\renewcommand{\arraystretch}{1}
\end{equation}
First it follows from equations (48) that the differential of an 
arbitrary point $x = x^0 A_0 + x^a A_a + x^{n+1} A_{n+1}$ of the 
generator $\alpha^{m+1} \subset U^n$ belongs to the tangent 
hyperplane $\xi$ of the hypersurface $U^n$  defined by the 
equation $x^n = 0$, and this confirms one more time the result 
of Theorem 5. 

Note  that  equations (48) allow us to find the 
equations of the focus surface ${\cal F}$ of the generator 
$\alpha^{m+1}$. In fact,  a {\em singular point} $x$ of the  
generator $\alpha^{m+1}$ is a point at which the dimension of the 
tangent subspace $T_x (U^n)$ is less than $n$. By (48), we can 
write this analytically as follows: 
\begin{equation}\label{eq:49}
dx \equiv (x^0 \omega_0^p  + x^a \omega_a^p  + x^{n+1} 
\omega_{n+1}^p) A_p \pmod{\alpha^{m+1}}.
\end{equation}
Since the forms $\omega_a^p$ are expressed in terms of the basis 
forms by formulas (45) and since by (5) and (46) the forms 
$\omega_{n+1}^p$ are expressed by the formulas 
$$
\omega_{n+1}^p = g^{pq}   \omega_q^0 = g^{pt} c_{tq} \omega^q 
= c_q^p \omega^q,
$$
where $c_q^p = g^{pt} c_{tq}$,  relation (49) takes the form
\begin{equation}\label{eq:50}
dx \equiv (x^0 \delta_q^p + x^a \lambda_{aq}^p  + x^{n+1} c_q^p) 
\omega^q A_p \pmod{\alpha^{m+1}}.
\end{equation}

The dimension of the tangent subspace to the hypersurface $U^n$ 
is lowered at those points in which the determinant of the matrix 
of coefficients of  relation (50) vanishes. Therefore the 
equation of the focus surface ${\cal F}$ in the generator 
$\alpha^{m+1}$ has the form 
\begin{equation}\label{eq:51}
\det (x^0 \delta_q^p + x^a \lambda_{aq}^p + x^{n+1} c_q^p) = 0,
\end{equation}
and is an algebraic equation of degree $r = n - m - 1$. Thus the 
hypersurface ${\cal F} \subset \alpha^{m+1}$ is an algebraic 
hypersurface of order $r$. 

The plane $\alpha^{m+1}$ intersects the hyperquadric $Q^n$ at the image of an $m$-dimensional sphere $S^m \subset C^n$ (a 
generator of the $m$-canal hypersurface $V^{n-1} \subset C^n$). 
By (2), the equation of this generator in $\alpha^{m+1}$ has the 
form 
\begin{equation}\label{eq:52}
(x, x) = g_{ab} x^a x^b - 2 x^0 x^{n+1} = 0.
\end{equation}
The intersection ${\cal F} \cap S^m$ defines singular points on 
the generator $S^m$ of the hypersurface $V^{n-1}$. 
Since we consider 
only the real geometry of canal hypersurfaces, we are interested 
in real singular points only. We denote the set of such points by 
$\Sigma$, $\Sigma = \mbox{{\rm Re}} ({\cal F} \cap S^m)$.

The classification of canal hypersurfaces $V^{n-1}$ is connected 
with the structure of the sets $\Sigma$ on the generators 
$S^m$. In general, such a set is an algebraic 
submanifold in $\alpha^{m+1}$ which has codimension 2 and order 
$2(n-m-1)$. However, in real geometry this set can be the empty 
set; then the generator  $S^m \subset V^{n-1}$ does not have 
singularities.

A detailed investigation of the structure of the set $\Sigma$ 
defined by equations (51) and (52) is a part of real 
algebraic geometry and falls outside of the scope of this paper. 

As an example, we consider a two-dimensional canal surface $V^2$ 
in the conformal space $C^3$ (see [S 92]). In this case we have 
$n = 3, m = 1$, and $r = 1$. The formulas of Sections {\bf 3} and 
{\bf 5} are still valid for this case but we should set $p, q, t 
= 2$ in these formulas. In addition, we can assume that $g_{22} = 
1$. The formulas (41) take the form 
$$
\omega_1^1 = \omega_2^2 = 0, \;\; \omega^1_2 + \omega_1^2 = 0.
$$
Next it follows from (39) that 
$$
\omega_3^1 = 0, \;\; \omega^2_3 = - \lambda_2^2 \omega^2,
$$
 and formula (44) becomes
$$
\omega_3^0 = 0.
$$
Thus we have
$$
dA_3 = - \lambda_2^2 \omega^2 A_2
$$
since $(dA_3, dA_3) = (\lambda_2^2 \omega^2)^2 > 0$, and as a 
result, the point $A_3$ describes a spacelike curve $l$ in the 
space $M^3_1$. 

Finally, formulas (45) and (46) take the form
$$
\omega_1^3 = \lambda_{12}^2 \omega^2,  \;\; \omega_2^0 
= c_{22} \omega^2.
$$

A two-dimensional plane generator $\alpha^2$ of a tangentially 
degenerate surface $U^3 \subset P^4$ that corresponds to a canal 
surface $V^2 \subset C^3$ under the Darboux mapping is defined by 
the equations $x^2 = x^3 = 0$. This and equation (52) imply that 
the equation of the generator $S^1$ of the canal surface $V^2$ 
has the form 
\begin{equation}\label{eq:53}
(x^1)^2 - 2 x^0 x^3 = 0.
\end{equation}

The focus hypersurface ${\cal F}$ in the generator $\alpha^2$ 
becomes a straight line that, by (51), is defined in $\alpha^2$ 
by the equation
\begin{equation}\label{eq:54}
x^0 + \lambda_{12}^2 x^1 + c_{22} x^4 = 0.
\end{equation}

The intersection of the focus straight line (54) with the 
generator (53) defines singular points of the generator. 
Excluding $x^0$ from equations (53) and (54), we arrive at 
the quadratic equation
$$
(x^1)^2 + \lambda_{12}^2 x^1 x^4 + c_{22} (x^4)^2 = 0.
$$
The solution of this equation is
$$
\frac{x^4}{x^1} = -\frac{1}{2} \lambda_{12}^2 \pm \frac{1}{2} 
\sqrt{D},
$$
where $D = (\lambda_{12}^2 )^2 - 4 c_{22}$. It follows that for 
$D > 0$, the generator $S^1$ carries two real singular points; 
for $D = 0$, it carries one double singular point; and 
for $D < 0$, it carries no real singular points.

Note also that if a space curve $l$,  described by the point 
$A_3 \subset P^4$ lying outside of the Darboux hyperquadric 
$Q^3$, 
is a closed curve, then a canal surface $V^2 \subset C^3$, whose 
image in $P^4$ is the curve $l$, is a tube (see [G 90]).

Consider, for example, a plane $\alpha^2$ in the space $P^4$. If 
this plane meets the hyperquadric $Q^3$ along a circle $\sigma$, 
then any spacelike closed curve $l$ lying in this plane defines 
in the space $C^3$ a smooth tube without selfintersections. In 
this case the curve $l$ is homotopic to the circle $\sigma$. If 
the plane $\alpha$ has no common points with the hyperquadric 
$Q^3$, then any smooth closed spacelike curve $l$ lying in this 
plane is spacelike and defines in the space $C^3$ a smooth tube, 
but this tube has selfintersections. Finally, if the plane 
$\alpha$ is tangent to  the hyperquadric $Q^3$ at a point $z$, 
then any smooth closed spacelike curve $l$  lying in this plane  
defines in the space $C^3$ a canal surface with the unique 
singular point $z$. In this case the point $z$ is located 
inside of the area bounded by  the curve $l$.

\end{document}